\documentclass[11pt]{amsart}

\usepackage{amsmath,amsthm,amsfonts,amssymb}

\usepackage{color}

\newtheorem{thm}[equation]{Theorem}
\newtheorem{cor}[equation]{Corollary}
\newtheorem{lem}[equation]{Lemma}
\newtheorem{prop}[equation]{Proposition}
\theoremstyle{definition}

\numberwithin{equation}{section}

\newcommand{\ann}{\hbox{\rm Ann}}
\definecolor{mjo}{rgb}{0,0,.9}

\newcommand{\cc}{\mathbb{C}}

\newcommand{\z}{\mathbb{Z}}

\newcommand{\setof}[2]{\{ #1 \, | \, #2 \} }

\newcommand{\witt}{\mathfrak W}

\newcommand{\n}{\mathfrak{n}}

\newcommand{\h}{\mathfrak{h}}

\newcommand{\B}{\mathfrak{b}}

\newcommand{\g}{\mathfrak{g}}

\newcommand{\p}{\mathcal{P}}

\newcommand{\pp}{\tilde{\mathcal{P}}}

\newcommand{\V}{\mathcal{V}}

\begin{document}

\title{Whittaker Modules for the Virasoro Algebra}

\author{Matthew Ondrus}
\address{Mathematics Department \\ Weber State University \\ Ogden, Utah  84408-1702}
\email{mattondrus@weber.edu}
\author{Emilie Wiesner}
\address{Department of Mathematics \\
           Ithaca College\\ Ithaca, New York 14850  }
\email{ewiesner@ithaca.edu}


\maketitle

\begin{abstract}
Whittaker modules have been well studied in the setting of complex semisimple Lie algebras.  Their definition can easily be generalized to certain other Lie algebras with triangular decomposition, including the Virasoro algebra. We define Whittaker modules for the Virasoro algebra and obtain analogues to several results from the classical setting, including a classification of simple Whittaker modules by central characters and composition series for general Whittaker modules.
\end{abstract}

\section{Introduction}
In this paper we define and investigate Whittaker modules for the Virasoro algebra.  The Virasoro algebra has been widely studied due in part to its interesting representation theory and its role in mathematical physics.  Whittaker modules were first studied in the setting of complex semisimple finite-dimensional Lie algebras.  Their definition can easily be generalized to certain other algebras with triangular decomposition, including the infinite-dimensional Virasoro algebra.

Whittaker modules were first discovered for $sl_2(\cc)$ by Arnal and Pinzcon \cite{AP}.  In \cite{Ko78}, Kostant defined Whittaker modules for all  finite-dimensional complex semisimple Lie algebras $\g$; his definition was motivated by the connection between these modules and the Whittaker equations in number theory.  In his paper,  Kostant classified the Whittaker modules of $\g$, demonstrating a strong connection with the center of $U(\g)$ among other results on Whittaker modules.  In \cite{Bl81}, Block showed that the simple modules for $sl_2(\cc)$ consist of highest (lowest) weight modules, Whittaker modules, and a third family obtained by localization.  Connections between Whittaker modules and Category ${\mathcal{O}}$  have also been made (see \cite{MS97} and \cite{Ba97}).

The definition of Whittaker modules is closely tied to the triangular decomposition of a finite-dimensional complex semisimple Lie algebra $\g$.  Because of this, it is natural to consider Whittaker modules for other algebras with a triangular decomposition.  Results for the complex semisimple Lie algebras have been extended to quantum groups by Sevoystanov, for $U_h(\g)$ (see \cite{Se00}), and Ondrus, for $U_q(sl_2)$ (see \cite{On05}).  (The quantum Serre relations preclude nontrivial Whittaker modules for other rational quantum groups $U_q(\g)$.)  Whittaker modules have also been studied for generalized Weyl algebras (see \cite{BO08}) and in connection to non-twisted affine Lie algebras (see \cite{Ch08}).

The Virasoro algebra $\V$ is the universal central extension of the Witt algebra $\witt$; $\witt$ is the smallest Cartan-type Lie algebra over $\cc$ and can be realized as the Lie algebra of polynomial vector fields on the circle.  The Virasoro algebra is an infinite-dimensional Lie algebra with a triangular decomposition (see \cite{MP95}). It arises in the representation theory of affine Kac-Moody Lie algebras:  if $\dot{\mathfrak g}$ is a simple finite-dimensional complex Lie algebra, then $\V \subseteq {\rm der} \left( {\mathcal L} (\dot{\mathfrak g}) \oplus \cc K \right)$ (see \cite{KacBook90}).  The Virasoro algebra is also closely connected to vertex operator algebras (cf. \cite{LL04}) and, through this topic, to areas of mathematical physics (cf.  \cite{Mandelstam}).

The triangular decomposition of the Virasoro algebra makes it possible to define Whittaker modules in this setting.  In this paper, we define Whittaker modules for the Virasoro algebra and obtain analogues to several of Kostant's results in the classical setting.  In Section \ref{sec:preliminaries}, we construct a universal Whittaker module for $\V$. This module is used to construct a Whittaker module with a central character, and we show in Corollary \ref{cor:L_PsiXiSimple} and Corollary \ref{cor:SimpleL_PsiXi} that this construction yields all simple Whittaker modules.  In Proposition \ref{Prop:zScalarImpliesSimple}, we show that every Whittaker module possessing a central character is simple (which gives a classification of Whittaker modules by central characters), and we describe the annihilator in $U(\V)$ of a cyclic Whittaker vector in this case.  In general, the action of the center of $U(\V)$ completely determines the structure of a Whittaker module, and in Theorem \ref{thm:generalV}, we describe the decomposition of an arbitrary Whittaker module based on the action of the center of $U(\V)$.  In Section \ref{sec:witt}, we use the results obtained for the Virasoro algebra to examine Whittaker modules for the Witt algebra.

We note that the proofs in the classical and Virasoro settings differ in their use of the center of the universal enveloping algebra.  Kostant's work utilizes knowledge of the center of $U(\g)$ from the Harish-Chandra homomorphism in order to characterize Whittaker vectors.   An equivalent tool does not exist for the Virasoro algebra.  Moreover, the center of $U(\V)$ is small (consisting only of polynomials of the central element $z$), and Fe{\u\i}gin and Fuks \cite{FF83} demonstrated the existence of commuting operators in the completion of $U(\V)$.  This suggests that the center of $U(\V)$ would not be sufficient to characterize the Whittaker modules for $U(\V)$, unlike the classical setting.  Surprisingly,  the results for the Virasoro algebra do in fact closely parallel those for complex semisimple Lie algebras.  Moreover, the concrete nature of the Virasoro algebra allowed us to  investigate Whittaker modules more directly in this setting. 

\section{Preliminaries} \label{sec:preliminaries}
Let $\V$ denote the Virasoro Lie algebra.  Then $\V={\rm span}_\cc \{z, d_k \mid k \in \z \}$
with Lie bracket 
\begin{align*}
[d_k, d_j] &= (j-k)d_{k+j} + \delta_{j,-k} \frac{k^3-k}{12} z; \\
[z, d_k] &= 0.
\end{align*}
We will make use of the following subalgebras:
\begin{eqnarray*}
\n^+ &= & {\rm span}_\cc \{ d_1, d_2, \ldots \} \\
\n^- &= & {\rm span}_\cc \{ d_{-1}, d_{-2}, \ldots \} \\
\h&=& {\rm span}_\cc \{ z, d_0 \} \\
\B^- &=& \n^- \oplus \h \\
\B^+ &=& \h \oplus \n^+
\end{eqnarray*}
Let $S(z)$ represent the symmetric algebra generated by $z$, that is, polynomials in $z$.  Then $S(z)$ is evidently contained in $Z(\V)$, the center of the universal enveloping algebra $U(\V)$. 
\subsection{Partitions and Pseudopartitions}
The following notation for partitions and pseudopartitions will be used to describe bases for $U(\V)$ and for Whittaker modules.

We define a {\it partition} $\mu$ to be a non-decreasing sequence of positive integers $\mu=( 0 < \mu_1 \leq \mu_2 \leq \cdots \leq \mu_r)$.  A {\it pseudopartition} $\lambda$ is a non-decreasing sequence of non-negative integers 
\begin{equation} \label{eqn:pseudopart1}
\lambda=( 0 \leq \lambda_1 \leq \lambda_2 \leq \cdots \leq \lambda_s).
\end{equation}
Let $\p$ represent the set of partitions, and let $\pp$ denote the set of pseudopartitions.  Then $\p \subseteq \pp$.

We also introduce an alternative notation for partitions and pseudopartitions.  For $\lambda \in \pp$, write
\begin{equation}\label{eqn:pseudopart2}
\lambda = (0^{\lambda(0)}, 1^{\lambda(1)}, 2^{\lambda(2)}, \ldots),
\end{equation}
where $\lambda(k)$ is the number of times $k$ appears in the pseudopartition and $\lambda(k)=0$ for $k$ sufficiently large.   Then a pseudopartition $\lambda$ is a partition exactly when $\lambda(0)=0$.  
For $\lambda \in \pp$, define
\begin{eqnarray*}
| \lambda | &=& \lambda_1 +  \lambda_2 + \cdots + \lambda_s \quad \mbox{(the size of $\lambda$)}\\
\# (\lambda) &=& \lambda(0) + \lambda(1) + \cdots \quad \mbox{(the $\#$ of parts of $\lambda$)}.
\end{eqnarray*}
For $\lambda \in \pp$, define elements $d_\lambda, d_{-\lambda} \in U(\V)$ by 
\begin{eqnarray*}
d_{ \lambda} &=& d_{\lambda_1} d_{\lambda_2}\cdots d_{\lambda_s} = d_0^{\lambda(0)} d_1^{\lambda(1)} \cdots \\
d_{-\lambda} &=& d_{-\lambda_s} \cdots d_{-\lambda_2}  d_{-\lambda_1}= \cdots d_{-1}^{\lambda(1)} d_0^{\lambda(0)}.
\end{eqnarray*}
Define $\overline{0}= (0^0, 1^0, 2^0, \ldots)$, and write $d_{\overline{0}}=1 \in U(\V)$.  We will consider $\overline{0}$ to be an element of $\pp$ but not of $\p$.
For any $\lambda \in \pp$ and $p(z) \in S(z)$, $p(z) d_{- \lambda}  \in U(\V)_{-|\lambda|}$, where $U(\V)_{-|\lambda|}$ is the $-|\lambda|$-weight space of $U(\V)$ under the adjoint action.
In particular, if $\lambda \in \p$, then $d_{- \lambda} \in U(\n^-)_{-|\lambda|}$.

\subsection{Whittaker Modules}\label{subsec:WhittakerModules}
In the classical setting of a finite-dimensional complex semisimple $\g$, a Whittaker module is defined  in terms of an algebra homomorphism from the positive nilpotent subalgebra $\g^+$ to $\cc$ (see \cite{Ko78}).  This homomorphism is required to be \emph{nonsingular}, meaning that it takes nonzero values on the Chevalley generators of $\g^+$.

In the present setting, the elements $d_1, d_2 \in \n^+$ generate $\n^+$.  Thus we assume that $\psi: \n^+ \rightarrow \cc$ is a Lie algebra homomorphism such that $\psi(d_1), \psi(d_2)\neq 0$, and we write $\psi_i = \psi (d_i)$ for $i > 0$.  We retain this assumption for the remainder of the paper.  The commutator relation in the definition of $\V$ forces $\psi_i = 0$ for $i \ge 3$.  For a $\V$-module $V$, a vector $w \in V$ is a {\it Whittaker vector} if $xw = \psi(x) w$ for all $x \in \n^+$.  A $\V$-module $V$ is a {\it Whittaker module} if there is a Whittaker vector $w \in V$ which generates $V$.  

For a given $\psi: \n^+ \rightarrow \cc$, define $\cc_\psi$ to be the one-dimensional $\n^+$-module given by the action $x \alpha = \psi (x) \alpha$ for $x \in \n^+$ and $\alpha \in \cc$.  The \emph{universal} Whittaker module $M_\psi$ is given by 
$$
M_{\psi}= U(\V) \otimes_{U(\n^+)} \cc_{\psi}.
$$
We use the term universal to refer to the property in Lemma \ref{lem:universalProp}(i) below.  Let $w=1 \otimes 1 \in M_\psi$.  By the PBW theorem $U( \B^-)$ has a basis $\setof{z^t  d_{- \lambda}}{\lambda \in \pp, t \in \z_{\ge 0}}$.  Thus $M_{\psi}$ has a basis 
\begin{equation} \label{eqn:PBWbasis}
\{z^t  d_{- \lambda} w \mid \lambda \in \pp, t \in \z_{\geq 0} \}
\end{equation}
and $uw \neq 0$ whenever $0 \neq u \in U( \B^-)$.  Define the {\it degree}  of $z^t d_{-\lambda}w$ to be $| \lambda |$.   For any $0 \neq v \in M_{\psi}$, define $maxdeg(v)$ to be the maximum degree of any nonzero component of homogeneous degree, and define $maxdeg(0) = -\infty$.  Define $max_{d_0} (v)$ to be the maximum value $\lambda (0)$ for any term $z^t d_{-\lambda} w$ with nonzero coefficient.   

For $\xi \in \cc$, define 
$$L_{\psi, \xi} = M_\psi / (z - \xi)M_\psi,$$
and let $\overline{\ \cdot \ }: M_{\psi} \rightarrow L_{\psi, \xi}$ be the canonical homomorphism.

\begin{lem}  \label{lem:universalProp}
Fix $\psi$ and $M_\psi$ as above.
\begin{itemize}
\item[(i)] Let  $V$ be a Whittaker module of type $\psi$ with cyclic Whittaker vector $w_V$.  Then there is a surjective map $\varphi: M_{\psi} \rightarrow V$ taking $w=1 \otimes 1$ to $w_V$.
\item[(ii)] Let $M$ be a Whittaker module of type $\psi$ with cyclic Whittaker vector $w_M$, and suppose that for any Whittaker module $V$ of type $\psi$ with cyclic Whittaker vector $w_V$ there exists a surjective homomorphism $\theta: M \rightarrow V$ with $\theta (w_M) = w_V$.  Then $M \cong M_{\psi}$.
\end{itemize}
\end{lem}
\begin{proof}
 To prove (i), note that if $v \in M_\psi$, then there exists $u \in U(\V)$ such that $v = uw$.  Define $\varphi : M_\psi \to V$ by $\varphi (v) = uw_V$.  To show that $\varphi$ is well-defined, it suffices to show that if $uw = 0$, then $uw_V = 0$.

For $u \in U(\V) = U( \B^-) \otimes U(\n^+)$, write $u = \sum_\alpha b_\alpha n_\alpha$ with $b_\alpha \in U(\B^-)$ and $n_\alpha \in U(\n^+)$.  Then $uw = \sum_\alpha b_\alpha n_\alpha w = \sum_\alpha \psi (n_\alpha) b_\alpha w$.  The PBW theorem implies that if $\sum_\alpha \psi (n_\alpha) b_\alpha \in U( \B^-)$ annihilates $w$, it must be that $\sum_\alpha \psi (n_\alpha) b_\alpha = 0$.  Since $w_V \in V$ is a Whittaker vector, it follows that 
$$0 = \left( \sum_\alpha \psi (n_\alpha) b_\alpha \right) w_V = \sum_\alpha b_\alpha \psi (n_\alpha) w_V = \sum_\alpha b_\alpha n_\alpha w_V = uw_V$$
as desired.

The proof of (ii) follows from (i).
\end{proof}

\section{Whittaker vectors in $M_{\psi}$ and $L_{\psi, \xi}$}
In this section we characterize the Whittaker vectors in $M_{\psi}$ and $L_{\psi, \xi}$.   
Fix $\psi : \n^+ \to \cc$ with $\psi_1, \psi_2 \neq 0$, and let $w = 1 \otimes 1 \in M_\psi$ as before.  For any $w' \in M_{\psi}$, $w'=u w$ for some $u \in U(\B^+)$, and thus $(d_n - \psi_n) w' = [d_n, u]w$.  Note that if $w'$ is a Whittaker vector, then $(d_n-\psi_n)w'=0$.

\begin{lem}\label{lem:computation2}
Define $M_\psi$ and $w = 1 \otimes 1 \in M_\psi$ as above, and let $a  \in \z_{>0}$, $k  \in \z_{\geq 0}$.  
Then, 
$$[d_{k+2}, d_{-k}^a]w = v-a (2k+2)\psi_2 d_{-k}^{a-1} w, $$
such that $maxdeg(v) <k(a-1)$ if $k>0$, and $max_{d_0} (v)<a-1$ if $k=0$.
\end{lem}
\begin{proof}
The lemma follows from a direct computation.
\end{proof}

\begin{lem}  \label{Lemma:deg}
Let $\overline{0} \neq \lambda \in \pp$.
\begin{itemize}
\item[(i)] For $m \in \z_{>0}$,  $maxdeg ([d_m, d_{-\lambda}] w) \leq |\lambda | -m +2$.
\item[(ii)] Suppose $k \in \z_{\geq 0}$ is minimal such that $\lambda(k) \neq 0$.  Then
$$
[d_{k+2}, d_{-\lambda}] w = v -  \lambda(k) \psi_2 (2k+2)\left(d_{-n}^{\lambda(n)} \cdots d_{-(k+1)}^{\lambda(k+1)} \right)d_{-k}^{\lambda(k)-1} w 
$$
where  if $k>0$, $maxdeg(v)<|\lambda|-k$; and if $k=0$, $v=v'+v''$ so that  $maxdeg(v')<|\lambda|-k$ and $max_{d_0} (v'')<\lambda(k)-1$.
\end{itemize}
\end{lem}
\begin{proof}
To prove (i), note that the commutator relations imply (using the notation from (\ref{eqn:pseudopart1})) that
\begin{eqnarray} \label{eqn:[dm,d_lambda]}
[d_m, d_{-\lambda}]  &=& \sum_i d_{-\lambda_n} \cdots [d_m, d_{-\lambda_i}]  \cdots d_{-\lambda_1} \nonumber \\
 &=&\sum_{0 \leq j \leq m} \ \sum_{\gamma \in \pp, |\gamma |= | \lambda | -m+j} p_{\gamma}(z) d_{-\gamma} d_j. \label{Equation:deg1}
\end{eqnarray}
The second equality in (\ref{eqn:[dm,d_lambda]}) results from the following argument.  If $m-\lambda_i  \leq 0$, then $[d_m, d_{-\lambda_i}] \in U(\B^-)$ and thus  $d_{-\lambda_n} \cdots [d_m, d_{-\lambda_i}]  \cdots d_{-\lambda_1} = \sum_{\gamma} p_\gamma (z) d_{-\gamma}$ for some $p_\gamma(z) \in S(z)$ and $\gamma \in \pp$ with $| \gamma| = |\lambda|-m$.  In the case that $m-\lambda_i >0$, we induct on $\# (\lambda)$ to obtain the stated form.   Since $d_j w=0$ for $j>2$, (i) follows from (\ref{Equation:deg1}).

To prove (ii) we first observe that (\ref{Equation:deg1}) can be refined to conclude that $j \leq m-\lambda_1$ for each $j$.  We may apply this reasoning to the partition 
$$\lambda'=(0^0, 1^0, \ldots, k^0, (k+1)^{\lambda(k+1)}, \ldots, m^{\lambda(m)}).$$
In the notation of (\ref{eqn:pseudopart1}), we have $\lambda'_1 = k+1$, and thus  
$$
[d_{k+2}, d_{-\lambda'}] d_{-k}^{\lambda(k)} w = \sum_{0 \leq j \leq 1} \sum_{\, \, \gamma \in \pp , \, |\gamma |= | \lambda' | -(k+2)+j} p_{\gamma} (z)d_{-\gamma} d_j d_{-k}^{\lambda(k)}w.
$$
 Therefore, 
 \begin{eqnarray}
 maxdeg ([d_{k+2}, d_{-\lambda'}] d_{-k}^{\lambda(k)} w) 
& \leq & max \{ |\gamma | \mid p_{\gamma}(z) \neq 0 \} +\lambda(k) \cdot k   \nonumber \\\
& \leq & |\lambda'|-(k+2)+1+\lambda(k) \cdot k \nonumber \\
& =&  |\lambda| -(k+1). \label{Equation:deg2}
\end{eqnarray}  
Note that
$[d_{k+2}, d_{-\lambda}] w =  [d_{k+2}, d_{-\lambda'}] d_{-k}^{\lambda(k)} w +d_{-\lambda'} [d_{k+2},  d_{-k}^{\lambda(k)}]  w$.  Statement (ii) follows from applying (\ref{Equation:deg2}) to the first summand and Lemma \ref{lem:computation2} to the second.  Note that the cases $k=0$ and $k>0$ come from the application of Lemma 3.1.
\end{proof}

\begin{prop}\label{prop:describeWhitVectsM}
Let $M_{\psi}$ be a universal Whittaker module for $\V$, generated by the Whittaker vector $w= 1 \otimes 1 $.  If $w' \in M_{\psi}$ is a Whittaker vector, then $w'=p(z)w$ for some $p(z) \in S(z)$.
\end{prop}
\begin{proof}  Let $w' \in M_{\psi}$ be an arbitrary vector.  By (\ref{eqn:PBWbasis}), 
$$w' = \sum_{\lambda \in \pp} p_{\lambda}(z)  d_{-\lambda} w$$
 for some polynomials $p_{\lambda}(z) \in S(z)$.   
 We will show that  if there is $\lambda \neq \overline{0}=(0^0, 1^0, 2^0, \ldots )$ such that $p_{\lambda}(z) \neq 0$, then there is $m \in \z_{>0}$ such that $(d_m - \psi_m)w'  \neq 0$.  In this case $w'$ is not a Whittaker vector, which proves the result.
   
Assume that $p_{\lambda}(z) \neq 0$ for at least one $\lambda \neq \overline{0}$.  Note that 
\begin{eqnarray}
(d_{m} - \psi_{m})w' 
&=& \sum_{\lambda \in \pp} [ d_{m}, p_{\lambda}(z)  d_{-\lambda}] w 
=  \sum_{\lambda \in \pp} p_{\lambda}(z) [ d_{m},  d_{-\lambda}] w, \label{Eq:w'}
\end{eqnarray}
where the second equality follows from the facts that $[d_{m}, -]$ is a derivation and $p_{\lambda}(z) \in Z(\V)$.  We shall argue that $(d_{m} - \psi_{m})w' \neq 0$ by considering each of the terms $[ d_{m},  d_{-\lambda}] w$.

Let $N= \mbox{max} \{ | \lambda | \mid  p_{\lambda}(z) \neq 0 \}$, and define $\Lambda_N = \{ \lambda \in \pp \mid  p_{\lambda}(z) \neq 0, |\lambda| = N \}$.  Assume $N>0$, and let $k$ be minimal such that $\mu (k) \neq 0$ for some $\mu \in \Lambda_N$.  It is clear that 
\begin{eqnarray}
(d_{k+2} - \psi_{k+2})w' 
&=& \sum_{\lambda \not\in \Lambda_N} p_{\lambda}(z) [ d_{k+2},  d_{-\lambda}] w  \label{Equation:prop35aa}\\
&& + \sum_{\lambda \in \Lambda_N} p_{\lambda}(z) [ d_{k+2},  d_{-\lambda}] w. \label{Equation:prop35ab}
\end{eqnarray}
We use this to deduce that 
\begin{eqnarray}
(d_{k+2} - \psi_{k+2})w' 
&=& v'+ v''  \label{Equation:prop35b} \\
& +& \sum_{\substack{\lambda \in \Lambda_N \\ \lambda(k) \neq 0}}   \lambda(k) (2k+2) \psi_2 p_{\lambda}(z) d_{-\lambda'}  w, \label{Equation:prop35c}
\end{eqnarray}
where $maxdeg (v')<N-k$; $max_{d_0} (v'') <\lambda(k)-1$; and $\lambda'$ is given by $\lambda'(k) = \lambda(k)-1$ and $\lambda'(i) = \lambda(i)$ for $i \neq k$.  To obtain this equality, we apply Lemma \ref{Lemma:deg} (i) to the sum in (\ref{Equation:prop35aa}) and Lemma \ref{Lemma:deg} (ii) to 
(\ref{Equation:prop35ab}). The terms in  (\ref{Equation:prop35c}) have degree $N-k$, and the power of $d_0$ is $\lambda(k)-1$.  Since the $d_{-\gamma}w$ form an $S(z)$-basis for $M_\psi$, this implies that the terms in  (\ref{Equation:prop35c}) are linearly independent from each other and from  (\ref{Equation:prop35b}).  Thus, $(d_{k+2} - \psi_{k+2})w' \neq 0$.
\end{proof}

\begin{prop} \label{Pr:SimpleWhittaker}
Let $w = 1 \otimes 1 \in M_{\psi}$ and $\overline{w}= \overline{1 \otimes 1} \in L_{\psi, \xi}$.  If $w' \in L_{\psi, \xi}$ is a Whittaker vector, then $w'=c\overline{w}$ for some $c \in \cc$.
\end{prop}

\begin{proof}
Note that the set $\{d_{- \lambda} \overline{w} \mid \lambda \in \pp \}$ spans $L_{\psi, \xi}$.  We claim that this set is linearly independent and thus a basis for $L_{\psi, \xi}$.  To check this, suppose there are  $a_{\lambda} \in \cc$, with at most finitely many $a_{\lambda} \neq 0$,  so that 
$$
0=\sum_{\lambda} a_{\lambda} d_{-\lambda} \overline{w}= \overline{ \sum_{\lambda} a_{\lambda} d_{-\lambda} w}
$$ 
in $L_{\psi, \xi}$.  In other words, $\sum_{\lambda} a_{\lambda} d_{-\lambda} w \in  (z - \xi)M_\psi$, and so
\begin{eqnarray*}
\sum_{\lambda} a_{\lambda} d_{-\lambda} w &=& (z - \xi) \sum_{\lambda, 0 \leq i \leq k} b_{\lambda, i} z^i d_{-\lambda} w 
\end{eqnarray*}
for some $k \in \z^+$ and $b_{\lambda,  i} \in \cc$.  This expression can be rewritten as 
$$
\sum_{\lambda} (a_{\lambda} - \xi b_{\lambda,0}) d_{-\lambda} w + \sum_{\lambda, 1 \leq i \leq k} (\xi b_{\lambda, i} - b_{\lambda, i-1} ) z^i d_{-\lambda} w - \sum_{\lambda} b_{\lambda, k} z^{k+1} d_{-\lambda}w=0
$$
Since each of the basis vectors (in $M_{\psi}$) in this linear combination are distinct, we conclude that $b_{\lambda, k}=0$, $\xi b_{\lambda, i} - b_{\lambda, i-1}=0$, $a_{\lambda} - \xi b_{\lambda,0}=0$, and thus $a_{\lambda}=0$ for all $\lambda \in \pp$.

With this fact now established, it is possible to use essentially the same argument as in Proposition \ref{prop:describeWhitVectsM} to complete the proof.  Here, however, we simply replace the polynomials $p_\lambda (z)$ in $z$ with scalars $p_\lambda$ whenever necessary.
\end{proof}

\section{Simple Whittaker Modules}
Using Proposition \ref{Pr:SimpleWhittaker}, we can now show that the modules $L_{\psi,\xi}$ are simple and form a complete set of simple Whittaker modules, up to isomorphism.

Fix an algebra homomorphism $\psi: \n^+ \rightarrow \cc$ and a Whittaker module $V$ of type $\psi$.
We may regard $V$ as an $\n^+$-module by restriction.  Define a modified action of $\n^+$ on $V$ (denoted by $\cdot$) by setting $x \cdot v = xv - \psi (x) v$ for $x \in \n^+$ and $v \in V$.  Thus if we regard a Whittaker module $V$ as an $\n^+$-module under the dot-action, it follows that $d_n \cdot v = d_nv - \psi_nv$ for $n>0$ and $v \in V$.

\begin{lem}\label{lem:locallyNilpotent}
  If $n>0$, then $d_n$ is locally nilpotent on $V$ under the dot action.
\end{lem}
\begin{proof}
Since $V = {\rm span}_{S(z)} \setof{d_{-\lambda}w}{\lambda \in \pp}$, it is sufficient to show that some power of $(d_n - \psi_n1)$ annihilates $d_{-\lambda}w$.  Note that $(d_n - \psi_n1)^k d_{-\lambda}w = {\rm ad}_{d_n}^k (d_{-\lambda})w$ and ${\rm ad}_{d_n}^k (d_{-\lambda}) \in U(\V)_{-|\lambda| + nk}$.  Also observe that ${\rm ad}_{d_n}^k (d_{-\lambda})$ is a sum of terms of the form $d_{a_1}d_{a_2} \cdots d_{a_s}$, where $s \le \# (\lambda)$ and $a_1 \le a_2 \le \cdots \le a_s \in \z$.  Now $d_{a_1}d_{a_2} \cdots d_{a_s}w = 0$ unless $a_s \le 2$.  However, if $a_1 \le \cdots \le a_s \le 2$, then $d_{a_1}d_{a_2} \cdots d_{a_s} \in U(\V)_m$ where $m = a_1 + \cdots + a_s \le 2s$.  Therefore by choosing $k$ sufficiently large so that $2 \cdot \# (\lambda) < -|\lambda| + nk$, we ensure that $a_s >2$, and thus $d_{a_1}d_{a_2} \cdots d_{a_s}w = 0$.
\end{proof}

\begin{lem}\label{lem:predotActionFinite}
Let $\lambda \in \pp$ and $i>0$.  
\begin{itemize}
\item[(i)] For all $n>0$, $d_n \cdot (z^i d_{-\lambda} w) \in \mbox{span}_{\cc} \{ z^j d_{-\mu} w \mid | \mu| + \mu(0) \leq |\lambda| + \lambda(0), j=i,i+1 \}$.
\item[(ii)]  If $n > |\lambda|+2$, then $d_n \cdot (z^i d_{-\lambda}w) = 0$.
\end{itemize}
\end{lem}
\begin{proof}
First note that $d_n \cdot (z^i d_{-\lambda} w)=z^i d_n \cdot ( d_{-\lambda} w)$.  Therefore, without loss of generality we assume that $i=0$.  

The proof of (i) is by induction on $|\lambda| + \lambda (0)$, with the case $|\lambda| + \lambda (0) = 0$ being obvious.  Moreover, 
if $| \lambda |=0$ the result follows from Lemma \ref{Lemma:deg},
so that it is no loss to assume that $|\lambda|>0$. Let $m = \max \setof{i}{\lambda (i)>0}$, and define $\lambda' \in \pp$ so that $d_{-\lambda} = d_{-m} d_{-\lambda'}$.    Then 
\begin{align}
d_n \cdot (d_{-\lambda} w) &= [d_n, d_{-\lambda}]w = [d_n, d_{-m}]d_{-\lambda'}w + d_{-m}[d_n, d_{-\lambda'}]w   \nonumber \\
&= -(m+n)d_{n-m}d_{-\lambda'}w + \delta_{m,n} \frac{n^3-n}{12}zd_{-\lambda'}w + d_{-m}[d_n, d_{-\lambda'}]w. \nonumber 
\end{align}
Since $|\lambda'| + \lambda' (0) \le |\lambda| + \lambda (0)$, it is clear that $zd_{-\lambda'}w$ has the desired form.  By induction, $[d_n,d_{-\lambda'}]w \in \mbox{span}_{\cc} \{ z^j d_{-\mu} w \mid | \mu| + \mu(0) \leq |\lambda'| + \lambda'(0), j =0,1 \}$.  It then follows that $d_{-m}[d_n, d_{-\lambda'}]w \in \mbox{span}_{\cc} \{ z^j d_{-\mu} w \mid | \mu| + \mu(0) \leq |\lambda| + \lambda(0), j=0,1 \}$ because $-m<0$, $|\lambda'| = |\lambda| - m$, and $\lambda'(0) = \lambda (0)$.  

To see why $d_{n-m}d_{-\lambda'}w$ has the stated form, we consider two cases.  If $n-m<0$, then straightening $d_{n-m}d_{-\lambda'}$ (i.e., rewriting  it in terms of the PBW basis of $U(\V)$) results in a linear combination of terms of the form $d_{-\mu}w$, where $|\mu| = m-n+|\lambda'| = -n+ |\lambda|$ and $\mu (0) = \lambda' (0) = \lambda (0)$.  If $n-m>0$, note that 
\begin{align*}
d_{n-m}d_{-\lambda'}w &= \psi_{n-m} d_{-\lambda'} w + d_{n-m} \cdot (d_{-\lambda'}w).
\end{align*}
It follows from induction that $d_{n-m} \cdot (d_{-\lambda'}w)$, and therefore $d_{n-m}d_{-\lambda'}w$, have the correct form.

The proof of  (ii) follows immediately from Lemma \ref{Lemma:deg} (i).
\end{proof}

\begin{lem}\label{lem:dotActionFinite}
Suppose $V$ is a Whittaker module for $\V$, and let $v \in V$.  Regarding $V$ as an $\n^+$-module under the dot action, $U(\n^+) \cdot v$ is a finite-dimensional submodule of $V$.
\end{lem}
\begin{proof}
Write $v = \sum_{i \geq 0, \lambda \in \pp} a_{i,\lambda} z^i d_{-\lambda}w$ for some $a_{i,\lambda} \in \cc$, and choose $k, \gamma$ so that $ | \gamma|+\gamma(0)+k$ is maximal with $a_{k, \gamma} \neq 0$.  For $n>0$, Lemma \ref{lem:predotActionFinite} implies that 
$$d_n \cdot v  \in \mbox{span}_{\cc} \{ z^j d_{-\mu} w \mid | \mu| + \mu(0)+j \leq |\gamma| + \gamma(0)+k+1 \}.$$ 
Since there are only finitely many pairs $(j, \mu)$ satisfying this condition, it follows that $U(\n^+) \cdot v$ is finite-dimensional.
\end{proof}

\begin{thm} \label{thm:submodule}
Let $V$ be a Whittaker module for $\V$, and let $S \subseteq V$ be a submodule.  Then there is a nonzero Whittaker vector $w' \in S$.
\end{thm}
\begin{proof}
Regard $V$ as an $\n^+$-module under the dot-action.  Let $0 \neq v \in S$, and let $F$ be the submodule of $S$ generated by $v$ under the dot-action of $\n^+$.  By Lemma  \ref{lem:dotActionFinite}, $F$ is a finite-dimensional $\n^+$-module.  Since Lemma \ref{lem:predotActionFinite} implies that $d_n \cdot F = 0$ for sufficiently large $n$, the quotient of $\n^+$ by the kernel of this action is also finite-dimensional.  Note that $d_n$ is locally nilpotent on $V$ (and thus on $F$) under this action by Lemma \ref{lem:locallyNilpotent}.  Thus Engel's Theorem implies that there exists a nonzero $w' \in F \subseteq S$ such that $x \cdot w'=0$ for all $x \in \n^+$.  By definition of the dot-action, $w'$ is a Whittaker vector.
\end{proof}

\begin{cor} \label{cor:L_PsiXiSimple}
For any $\xi \in \cc$, $L_{\psi, \xi}$ is simple.  
\end{cor}
\begin{proof}
Let $S$ be a nonzero submodule of $L_{\psi, \xi}$.  Since $z \in \V$ acts by the scalar $\xi$ on $L_{\psi, \xi}$, it follows from Theorem \ref{thm:submodule} that there is a nonzero Whittaker vector $w' \in S$.  Proposition \ref{Pr:SimpleWhittaker} implies that $w' = cw$ for some $c \in \cc$, and therefore $w \in S$.  Since $w$ generates $L_{\psi, \xi}$, we have $S= L_{\psi, \xi}$.
\end{proof}

\noindent 
\begin{cor} \label{cor:SimpleL_PsiXi}
Let $\psi : \n^+ \to \cc$ be a Lie algebra homomorphism such that $\psi (d_1), \psi (d_2) \neq 0$, and let $S$ be a simple Whittaker module of type $\psi$ for $\V$. Then $S \cong L_{\psi, \xi}$ for some $\xi \in \cc$. 
\end{cor}
\begin{proof}
Let $w_s \in S$ be a cyclic Whittaker vector corresponding to $\psi$.  By Schur's lemma, the center of $U( \V )$ acts by a scalar, and this implies that there exists $\xi \in \cc$ such that $z s = \xi  s$ for all $s \in S$.  Now by the universal property of $M_\psi$, there exists a module homomorphism $\varphi : M_\psi \to S$ with $uw \mapsto uw_s$.  This map is surjective since $w_s$ generates $S$.  But then $\varphi \left( (z - \xi) M_\psi \right) = (z - \xi) \varphi (M_\psi ) = (z - \xi)S = 0$, so it follows that 
$$(z - \xi) M_\psi \subseteq \ker \varphi \subseteq M_\psi.$$
Because $L_{\psi, \xi}$ is simple and $\ker \varphi \neq M_\psi$, this forces $\ker \varphi = (z - \xi) M_\psi$.
\end{proof}

For a given $\psi : \n^+ \to \cc$ and $\xi \in \cc$, note that $L = U(\V) (z-\xi1) + \sum_{i>0} U(\V)(d_i - \psi_i1) \subseteq U(\V)$ is a left ideal of $U(\V)$.  For $u \in U(\V)$, let $\overline u$ denote the coset $u + L \in U(\V)/L$.  Then we may regard $U(\V)/L$ as a Whittaker module of type $\psi$ with cyclic Whittaker vector $\overline 1$.  

\begin{lem}\label{lem:quotientEnveloping}
Fix $\psi : \n^+ \to \cc$ with $\psi_1, \psi_2 \neq 0$.  Define the left ideal $L$ of $U(\V)$ by $L = U(\V) (z-\xi1) + \sum_{i>0} U(\V)(d_i - \psi_i1)$, and regard $V = U(\V) / L$ as a left $U(\V)$-module.  Then $V$ is simple, and thus $V \cong L_{\psi, \xi}$.  
\end{lem}
\begin{proof}
Note that the center of $U( \V )$ acts by the scalar $\xi$ on $V$.  By the universal property of $M_\psi$, there exists a module homomorphism $\varphi : M_\psi \to V$ with $uw \mapsto u\overline 1$.  This map is surjective since $\overline 1$ generates $V$.  But then $\varphi \left( (z - \xi) M_\psi \right) = (z - \xi) \varphi (M_\psi ) = (z - \xi)V = 0$, so it follows that 
$$(z - \xi) M_\psi \subseteq \ker \varphi \subseteq M_\psi.$$
Because $L_{\psi, \xi}$ is simple and $\ker \varphi \neq M_\psi$, this forces $\ker \varphi = (z - \xi) M_\psi$.
\end{proof}

\begin{prop}\label{Prop:zScalarImpliesSimple}
Suppose that $V$ is a Whittaker module of type $\psi$ such that $z \in \V$ acts by the scalar $\xi \in \cc$.   Then $V$ is simple.  Moreover, if $w$ is a cyclic Whittaker vector for $V$, $\ann_{U(\V)} (w) = U(\V) (z-\xi1) + \sum_{i>0} U(\V)(d_i - \psi_i1)$.
\end{prop}
\begin{proof}
Let $K$ denote the kernel of the natural surjective map $U(\V) \to V$ given by $u \mapsto uw$.  Then $K$ is a proper left ideal containing $L = U(\V) (z-\xi1) + \sum_{i>0} U(\V)(d_i - \psi_i1)$.  By Lemma \ref{lem:quotientEnveloping}, $L$ is maximal, and thus $K=L$ and $V \cong U(\V)/L$ is simple.
\end{proof}

Note that by Schur's Lemma, Proposition \ref{Prop:zScalarImpliesSimple} applies to any simple Whittaker module.

\section{Arbitrary Whittaker modules}
We now characterize arbitrary Whittaker modules, with generating Whittaker vector $w$, in terms of the annihilator $\ann_{S(z)} (w)$.

\begin{lem}\label{lem:uniserial}
Suppose that $V$ is a Whittaker module of type $\psi$ with cyclic Whittaker vector $w$, and assume that $\ann_{S(z)}(w) = (z- \xi1)^a$ for some $a>0$.  Define submodules $V = V_0 \supseteq V_1 \supseteq \cdots \supseteq V_a = 0$ by $V_i = U(\V)(z- \xi1)^iw$.  Then 
\begin{itemize}
  \item[(i)] $V_i$ is a Whittaker module of type $\psi$, with cyclic Whittaker vector $w_i = (z- \xi1)^iw$, 
  \item[(ii)] $V_i/V_{i+1}$ is simple for $0 \le i < a$, and 
  \item[(iii)] the submodules $V_0, \ldots, V_a$ are the only submodules of $V$. 
\end{itemize}
\end{lem}
\begin{proof}
It is immediate that $V_i$ is a Whittaker module with cyclic Whittaker vector $w_i$.  Since $z$ acts by the scalar $\xi$ on $V_i/V_{i+1}$, it follows from Proposition \ref{Prop:zScalarImpliesSimple} that the $V_i/V_{i+1}$ is simple, and thus isomorphic to $L_{\psi, \xi}$.  Therefore $V_0 \supseteq V_1 \supseteq \cdots \supseteq V_a$ form a composition series for $V$, and any simple subquotient of $V$ must be isomorphic to $V/V_1 \cong L_{\psi, \xi}$.

If $M$ is any maximal submodule of $V$, then $V/M$ is simple and thus $z$ acts by some scalar $\kappa \in \cc$. On the other hand, $(z - \xi 1)^a$ acts as 0 on $V$, and therefore on $V/M$, so that $\kappa = \xi$.   This implies that $(z- \xi1)V \subseteq M$.  Since $V_1 = U(\V)(z- \xi1)w = (z- \xi1)V \subseteq M$, it follows from the maximality of $V_1$ and $M$ that $V_1=M$.  A similar argument shows that $V_{i+1}$ is the unique maximal submodule of $V_i$ for every $i<a$, and thus the $V_i$ are the only submodules of $V$.
\end{proof}

\begin{thm} \label{thm:generalV}
Assume that $V$ is a Whittaker module of type $\psi$ with cyclic Whittaker vector $w$, and suppose that $\ann_{S(z)}(w) \neq 0$.  

Let $p(z)$ be the unique monic generator of the ideal $\ann_{S(z)}(w)$ in $S(z)$, and write $p(z) = \prod_{i=1}^k (z- \xi_i1)^{a_i}$ for distinct $\xi_1, \ldots, \xi_k \in \cc$.  Define 
$$
w_j = p_j(z)w, \ \mbox{where} \  p_j(z) = \prod_{i \neq j} (z- \xi_i1)^{a_i}, \quad \mbox{and} \ V_j = U(\V)w_j.
$$ 
Then $V_i$ is a Whittaker module of type $\psi$ with cyclic Whittaker vector $w_i$, and $V = V_1 \oplus \cdots \oplus V_k$.  Furthermore, the submodules $V_1, \ldots, V_k$ are indecomposable; $V_j$ is simple if and only if $a_j = 1$; and $a_j$ is the composition length of $V_j$.  
\end{thm}
\begin{proof}
Since $\gcd (p_1(z), \ldots, p_k(z)) = 1$, there exist polynomials $q_1(z), \ldots, q_k(z) \in S(z)$ such that $\sum_i q_i(z)p_i(z) = 1$.  This implies that $1w = \left( \sum_i q_i(z)p_i(z) \right) w \in V_1 + \cdots + V_k$, and thus $V = V_1 + \cdots + V_k$.

To show that the sum $V_1 + \cdots + V_k$ is direct, first note that for $i \neq j$, $p(z)$ divides $p_i(z)p_j(z)$, and thus $p_j(z)w_i = 0$.  It follows that 
\begin{align*}
w_i &= 1w_i \\
&= (q_1(z)p_1(z) + \cdots + q_k(z)p_k(z))w_i \\
&= q_i(z)p_i(z)w_i.
\end{align*}
Now if there exist $u_1, \ldots, u_k \in U(\V)$ such that $u_1w_1+ \cdots + u_kw_k = 0$, then 
$$0 = q_i(z)p_i(z) \left( \sum_j u_jw_j \right) = u_i q_i(z)p_i(z)w_i = u_iw_i,$$
and this implies that the sum is direct.

That the submodules $V_1, \ldots, V_k$ are indecomposable with the stated composition length follows from Lemma \ref{lem:uniserial}.
\end{proof}

\begin{cor}
Let $V$ be a Whittaker module of type $\psi$ with cyclic Whittaker vector $w$, assume $\ann_{S(z)}(w) \neq 0$.  Let $p(z)$ be the unique monic generator of $\ann_{S(z)}(w)$.  Then $\ann_{U(\V)} (w) = U(\V) p(z) + \sum_{i>0} U(\V)(d_i - \psi_i1)$.
\end{cor}
\begin{proof}
We use induction on the composition length of $V$.  Write $p(z) = (z-\xi1)p'(z)$ for some $\xi \in \cc$ and some monic polynomial $p'(z) \in S(z)$.  If $p'(z) = 1$, then $V$ is simple, and the result is true by Proposition \ref{Prop:zScalarImpliesSimple}.  Thus we assume that $(z-\xi1)p'(z)$ is a nontrivial factorization of $p(z)$ and $(z-\xi1)w \neq 0$.  

Let $w' = (z-\xi1)w$, and let $V' = U(\V)w' \subseteq V$.  Then $V'$ is a Whittaker module with cyclic Whittaker vector $w'$, and $\ann_{S(z)}(w') = S(z)p'(z)$.  Theorem \ref{thm:generalV} therefore implies that the composition length of $V'$ is one less than that of $V$, and it follows by induction that $\ann_{U(\V)}(w') = U(\V) p'(z) + \sum_{i>0} U(\V)(d_i - \psi_i1)$.  Let $\overline w = w + V' \in V/V'$, and observe that $\ann_{S(z)}(\overline w) = S(z) (z-\xi1)$.

Let $u \in \ann_{U(\V)}(w)$.  Since $\ann_{U(\V)}(w) \subseteq \ann_{U(\V)}(\overline w)$, Proposition \ref{Prop:zScalarImpliesSimple} implies that 
\begin{equation} \label{eqn:uAnnhilatewBar}
u = u_0(z-\xi1) + \sum_{i \ge 1} u_i(d -\psi_i1) \in U(\V) (z-\xi1) + \sum_{i>0} U(\V)(d_i - \psi_i1)
\end{equation}
But $\sum_{i \ge 1} u_i(d -\psi_i1) \in \ann_{U(\V)}(w)$, so $u_0(z-\xi1) \in \ann_{U(\V)}(w)$.  Observe that $0 = u_0(z-\xi1)w = u_0w'$, and thus 
$$u_0 \in \ann_{U(\V)}(w') = U(\V)p'(z) + \sum_{i>0} U(\V)(d_i - \psi_i1).$$
It follows from (\ref{eqn:uAnnhilatewBar}) that $u$ has the required form.
\end{proof}

\begin{lem} \label{lem:AnnZero}
Let $V$ be a Whittaker module for $\V$ with cyclic Whittaker vector $w$, and assume that $\ann_{S(z)} (w) = 0$.  Then $V \cong M_\psi$.
\end{lem}
\begin{proof}
There is a surjective homomorphism $\varphi: M_{\psi} \rightarrow V$.  The kernel of this map, $K$, is a submodule of $M_{\psi}$.  If $K \neq 0$, then by Theorem \ref{thm:submodule}, we know that there is a nonzero Whittaker vector $w' \in K$.  Proposition \ref{prop:describeWhitVectsM} implies that $0 \neq w'=p(z) 1 \otimes 1$ and thus $0 \neq p(z) \in \ann_{S(z)} (w)$, which is a contradiction.  Therefore, it must be that $K = 0$ and $\varphi$ is an isomorphism.
\end{proof}

\begin{thm}
Let $M_\psi$ be the universal Whittaker module of type $\psi$, and let $w = 1 \otimes 1 \in M_\psi$.  If $V \subseteq M_\psi$ is a submodule, then $V \cong M_\psi$.  Furthermore, $V$ is generated by a Whittaker vector of the form $q(z)w$ for some polynomial $q(z)$.
\end{thm}
\begin{proof}
Recall that for $\xi \in \cc$, the submodule of $M_\psi$ generated by $(z - \xi)w$ is maximal.  Note that $U(\V)(z - \xi)w$ is a Whittaker module with cyclic Whittaker vector $(z - \xi)w$ (of type $\psi$), and by Lemma \ref{lem:AnnZero} $U(\V)(z - \xi)w \cong M_\psi$.  

Now by Proposition \ref{prop:describeWhitVectsM} and Theorem \ref{thm:submodule}, the given submodule $V$ contains a Whittaker vector of the form $p(z)w$, and by Lemma \ref{lem:AnnZero}, the submodule $M'$ generated by $p(z)w$ is isomorphic to $M_\psi$.  As $\cc$ is algebraically closed and $p(z)$ can be written as a product of linear factors, the first paragraph implies the existence of a chain of universal Whittaker modules between $M'$ and $M_\psi$ so that each quotient is irreducible.  This implies that $V$ must be one of the submodules in the chain, so in fact $V$ is a universal Whittaker module of type $\psi$.
\end{proof}

\begin{cor}
Let $V$ be a Whittaker module of type $\psi$, with cyclic Whittaker vector $w$.  Then ${\rm Wh}(V) = S(z) w$.
\end{cor}
\begin{proof}
 If $\ann_{S(z)} (w)=0$, the result follows from Proposition \ref{prop:describeWhitVectsM} and Lemma \ref{lem:AnnZero}.  Therefore, we assume that $\ann_{S(z)}(w) \neq 0$.  In this case, Theorem \ref{thm:generalV} implies that $V$ has finite composition length $a$.  We will prove the result by induction on $a$.

If $a=1$, then $V$ is a simple module, and the result follows from Proposition \ref{Pr:SimpleWhittaker}.  Now, suppose $V$ is a module with arbitrary composition length $a$.  Without loss of generality, we may assume that $V$ is indecomposable, so that the composition series is of the form  $V = V_0 \supseteq V_1 \supseteq \cdots \supseteq V_a = 0$ where $V_i$ has cyclic Whittaker vector $w_i = (z- \xi1)^iw$.  (See Lemma \ref{lem:uniserial} .)  

Let $w' \in V$ be a Whittaker vector.  Since $V/V_1$ is simple, Proposition  \ref{Pr:SimpleWhittaker} implies that the image of $w'$ in $V/V_1$ will be a scalar multiple of the image of $w$.  Therefore, in $V$, $w'=cw + w''$ for some $c \in \cc$ and $w'' \in V_1$.  Note that $w''=w'-cw$ is also a Whittaker vector.  Since $V_1$ has composition length $a-1$, it follows from the inductive hypothesis that $w'' = p(z)w_1=p(z) (z- \xi1)w$ for some $p(z) \in S(z)$.  Therefore, $w'=cw+p(z)(z- \xi1)w$, the desired result.
\end{proof}

\section{Whittaker modules for the Witt algebra} \label{sec:witt}

In this section, we describe the Whittaker modules for the Witt algebra $\witt$.  Recall that the Virasoro algebra $\V$ is the universal central extension of $\witt$.  We abuse notation and regard $\witt={\rm span}_\cc \{ d_k \mid k \in \z \}$
with Lie bracket given by 
\begin{align*}
[d_k, d_j] &= (j-k)d_{k+j}
\end{align*}
for $j, k \in \z$.  As $\V$ is the universal central extension of $\witt$, there is a surjective Lie algebra homomorphism $\rho : \V \to \witt$ with $\ker \rho=\cc z$.  This map extends to a surjective homomorphism $U(\V) \to U(\witt)$ which we also denote by $\rho$.

Define the subalgebra $\n^+_\witt \subseteq \witt$ in the obvious manner.  Since $\n^+ \cong \n^+_\witt$, we make no distinction between a homomorphism $\psi : \n^+_\witt \to \cc$ and a homomorphism $\psi : \n^+ \to \cc$.  Let $\psi: \n^+_\witt \rightarrow \cc$ be a Lie algebra homomorphism such that $\psi(d_1), \psi(d_2)\neq 0$.  A $\witt$-module $V$ is a {\it Whittaker module} if there is some $w \in V$ such that $w$ generates $V$ and $xw=\psi(x)w$ for all $x \in \n^+_\witt$.

\begin{prop}
Fix a homomorphism $\psi : \n^+_\witt \to \cc$ with $\psi (d_1), \psi (d_2) \neq 0$, and let $V$ be a nonzero Whittaker module of type $\psi$ for $\witt$.  Then $V$ is simple.  Moreover, $V \cong L_{\psi, 0}$ when $L_{\psi, 0}$ is viewed as a $\witt$-module.
\end{prop}
\begin{proof}
Let $V_\V$ be the $\V$-module obtained by letting $x \in \V$ act on $V$ by $\rho (x) \in \witt$.  Then $V_\V$ is a nonzero Whittaker module for $\V$, and the central element $z \in \V$ acts by $0$.  
By Proposition \ref{Prop:zScalarImpliesSimple}, we then have $V_{\V} \cong L_{\psi,0}$.  As $V_\V$ is the pullback of $V$ and $\rho : \V \to \witt$ is surjective, we conclude that $V$ must also be simple.  

To check that $L_{\psi, 0}$ can be viewed at a $\witt$-module, we note that $z \in \ann_{\V} (L_{\psi, 0})$.  Therefore, the action of $\witt = \V/ \cc z$ is well-defined.  Since $V_{\V} \cong L_{\psi,0}$ as $\V$ modules, we must have $V \cong L_{\psi,0}$ as $\witt$-modules.
\end{proof}

\end{document}